\documentclass[12pt]{article}
\usepackage{amsfonts}
\usepackage{latexsym}

\newcommand{\be}{\begin{equation}} \newcommand{\ee}{\end{equation}}
\newcommand{\bea}{\begin{eqnarray}} \newcommand{\eea}{\end{eqnarray}}
\newcommand{\bean}{\begin{eqnarray*}}
  \newcommand{\eean}{\end{eqnarray*}}
\newcommand{\brray}{\begin{array}} \newcommand{\erray}{\end{array}}
\newcommand{\ben}{\begin{equation}{nonumber}}
  \newcommand{\een}{\end{equation}{nonumber}}
\newcommand{\newsection}[1]{\setcounter{equation}{0}
  \setcounter{dfn}{0}
\section{#1}}

\newtheorem{dfn}{Definition}[section] \newtheorem{thm}[dfn]{Theorem}
\newtheorem{lmma}[dfn]{Lemma} \newtheorem{ppsn}[dfn]{Proposition}
\newtheorem{crlre}[dfn]{Corollary} \newtheorem{xmpl}[dfn]{Example}
\newtheorem{rmrk}[dfn]{Remark}

\newcommand{\bdfn}{\begin{dfn}} \newcommand{\bthm}{\begin{thm}}
    \newcommand{\blmma}{\begin{lmma}}
      \newcommand{\bppsn}{\begin{ppsn}}
        \newcommand{\bcrlre}{\begin{crlre}}
          \newcommand{\bxmpl}{\begin{xmpl}}
            \newcommand{\brmrk}{\begin{rmrk}}

              \newcommand{\edfn}{\end{dfn}}
            \newcommand{\ethm}{\end{thm}}
          \newcommand{\elmma}{\end{lmma}}
        \newcommand{\eppsn}{\end{ppsn}}
      \newcommand{\ecrlre}{\end{crlre}}
    \newcommand{\exmpl}{\end{xmpl}} \newcommand{\ermrk}{\end{rmrk}}

\newcommand{\IC}{{\Bbb C}}

 \newcommand{\IR}{{\Bbb R}}
 \newcommand{\IT}{{\Bbb T}}

 \newcommand{\IZ}{{\Bbb Z}}

\newcommand{\ath}{{\cal A}_\theta} \newcommand{\al}{\alpha}
\newcommand{\bta}{\beta} \newcommand{\gma}{\gamma}
\newcommand{\Dlt}{\Delta} \newcommand {\Nab}{\nabla}
\newcommand{\dlt}{\delta} 
 
\newcommand{\lmd}{\lambda}

 \newcommand{\tta}{\theta}

\newcommand{\cla}{{\cal A}} \newcommand{\clb}{{\cal B}}
 \newcommand{\cld}{{\cal D}}
\newcommand{\cle}{{\cal E}} \newcommand{\clf}{{\cal F}}
 \newcommand{\clh}{{\cal H}}
 
\newcommand{\clk}{{\cal K}} \newcommand{\cll}{{\cal L}}
 \newcommand{\cln}{{\cal N}}
 
 \newcommand{\cls}{{\cal S}}
\newcommand{\clt}{{\cal T}}

  \def 
\bbE
{\mbox{\boldmath $E$}}                     

  \def 
\bbe
{\mbox{\boldmath $e$}}

 \def\a*{{\cal A}_{h,*}} \def\B{{\cal B}(h)}
\def\B1{{\cal B}_1(h)} \def\b{{\cal B}^{s. a. }(h)} \def\b1{{\cal
    B}^{s. a. }_1(h)}

 \newcommand{\ot}{\otimes}
 
\newcommand{\raro}{\rightarrow}

 \newcommand {\CC}{\centerline} \def \qed 
{
  \mbox{}\hfill $\Box$\vspace{1ex}} 

\begin{document}

\CC {\large {\bf Probability and Geometry on some Noncommutative
    Manifolds } } \CC {by} \CC{Partha Sarathi Chakraborty,} \CC
{Debashish Goswami } \CC {and} \CC {Kalyan B. Sinha }
\vspace{.5 in}
{\it Abstract :} In a noncommutative torus, effect of perturbation by
inner derivation on the associated quantum stochastic process and
geometric parameters like volume and scalar curvature have been
studied. Cohomological calculations show that the above perturbation
produces new spectral triples.  Also for the Weyl $C^*$-algebra, the
Laplacian associated with a natural
stochastic process is obtained and associated volume form is 
calculated.
AMS classification (1991): 46L87, 81S25.
\newsection{Introduction} For a fixed $\theta$, an irrational number
in $[0,1]$, consider the $C^*$-algebra $\ath$ generated by a pair of
unitary symbols subject to the relation : \bea {\label 1}
\lefteqn{UV}~~ &=& {\rm exp}(2 \pi i \tta) VU \equiv \lmd VU.  \eea
For details of the properties of such a $C^*$-algebra, the reader is
referred to
\cite{Con}, \cite{Rf}. The algebra has many interesting 
representations :\\
(i) $\clh =L^2(\IT^1)$, $\IT^1$ is the circle, and for $f \in \clh ,$ 
$(\pi_1(U)f)(z)=f(\lmd z),$\\
$(\pi_1(V)f)(z)=zf(z),~z \in \IT^1.$\\
(ii) In the same $\clh$, with the roles of $U$ and $V$ reversed :\\
for $f \in \clh ,$ $(\pi_2(V)f)(z)=f(\bar{\lmd} z),
~(\pi_2(U)f)(z)=zf(z),~z \in \IT^1.$\\
(iii) In $\clh=L^2(\IR),$\\
$(\pi_3(U)f)(x)=f(x+1),~ (\pi_3(V)f)(x)=\lmd^x f(x).$\\
While the first two were inequivalent irreducible representations, 
the
ultra-weak closure of the third one is a factor of type $II_1$.

There is a natural action of the abelian compact group $\IT^2$
(2-torus) on $\ath$ given by,
$$\al_{(z_1,z_2)}(\sum a_{mn} U^mV^n)=\sum a_{mn}z_1^m z_2^n U^m
V^n,$$
where the sum is over finitely many terms and $\| z_1 \|=\|z_2
\|=1.$ $\al $ extends as a $*$-automorphism on $\ath$ and has two
commuting generators $d_1$ and $d_2$ which are skew-*-derivations
obtained by extending linearly the rule: \bea {\label 3}
{d_1(U)}= U,~d_1(V)=0 {\nonumber}\\
d_2(U)=0,d_2(V)=V.  \eea Both $d_1$ and $d_2$ are clearly well 
defined
on ${\ath}^\infty \equiv \{ a \in \ath \mid z \mapsto \al_z(a) $ is
$C^{\infty} \} \equiv \{ {\sum}_{m,n \in Z } a_{mn} U_m V^n \mid
sup_{m,n} | m^k n^l  a_{mn}| \leq c $ for all $k,l \in N \}.$ Since
the action is norm continuous $~\ath^{\infty}~ $ is a dense
$*$-subalgebra of $\ath$.  A theorem of Bratteli, Elliot, Jorgensen
\cite{BEJ} describes all the derivaions of $\ath $ which maps
$\ath^\infty $ to itself : for almost all $\theta$ (Lebesgue), a
derivation $\delta : \ath^\infty \rightarrow \ath^\infty $ is of the
form $\delta = c_1 d_1 +c_2 d_2 + [r,.]$, with $r \in \ath^\infty
,~c_1,c_2 \in \IC .$ Another important fact about $\ath$ is the
existence of a unique faithful trace $\tau$ on $\ath $ defined as
follows: \bea
 \label {4}
 \tau ( \sum a_{mn} U^m V^n )=a_{00}.  \eea Then one can consider the
 Hilbert space $ \clh = L^2 ( \ath , \tau ) $ (see \cite {Nel} for an
 account on noncommutative $L^p $ spaces.)  and study the derivations
 there. It is easy to see that the family $\{ U^mV^n \}_{m,n \in Z } 
$
 constitute a complete orthonormal basis in $\clh$. The next simple
 theorem is stated without proof.  \bthm The canonical derivations
 $d_1,d_2$ are self adjoint on their natural domains: $Dom(d_1)= \{
 \sum a_{mn} U^m V^n \mid \sum (1+  m^2 ) | a_{mn} |^2 < \infty \}$
 $Dom(d_2)= \{ \sum a_{mn} U^m V^n \mid \sum (1+
 n^2 ) | a_{mn}
 |^2 < \infty \}$. Furthermore if we denote by $d_r = [r,.] $ with $r
 \in \ath \subset L^{\infty} ( \ath , \tau ) $ acting as left
 multiplication in $\clh$, then $ {d_r}^*=d_{r^*} \in \clb ( \clh ) $
 \ethm

 \newsection{ Diffusion on $\ath $ and a noncommutative Laplacian }
 There is a canonical construction of a quantum stochastic flow or
 diffusion on a von Neumann \cite {GS} or a $C^*$-algebra $\cla $
 \cite{GSP} associated with a completely positive semigroup on $\cla
 $. The question about which of these semigroups have `local'
 generators $\cll $ remains open , though Sauvageot studied these in
 \cite {Sav} . Following these studies , we know that
 $\cll $  is characterized by  : \\
 (i) $\cld \subseteq Dom ( \cll ) \subseteq \cla \subseteq \clb ( 
\clh
 ) $, dense in
 $\cla $ such that  $\cld $ itself is a $*$-algebra,\\
 (ii) a $*$-representation $\pi$ in some Hilbert space $\clk $ and an
 associated $\pi $ derivation $\dlt $ such that $ \dlt ( x ) \in \clb
 ( \clh , \clk ) $ and $ \dlt
 (xy)= \dlt (x) y + \pi (x ) \dlt ( y ) $, \\
 (iii) a second order cocycle relation : $ \cll ( x^* y) - { \cll ( x
   ) }^* y -x^* \cll ( y ) = {\dlt (x ) }^* \dlt ( y ) $, for $ x,y
 \in \cld $.  In analogy with the heat semigroup in the case of
 classical diffusion, we shall call $\cll $
 the non-commutative Laplacian or  Lindbladian.  \\
 Hudson and Robinson \cite {HR} studied the above question for $\ath$
 in the case where the representation $\pi$ is the identity
 representation in $\clh $ itself and concluded that while there 
exist
 classical stochastic dilations for the Lindbladians $ \cll (x) = -
 \frac {1} {2} {d_1}^2 (x) $ or $ - \frac {1} {2}{d_2}^2 (x) $, there
 does not exist any $\cll $ corresponding to $\dlt = d_1 + i d_2 $ so
 that there is no quantum stochastic dilation corresponding to this
 case. We claim that if we choose $\pi (x ) = x \otimes I $ in $\clk 
=
 \clh \otimes C^2 \cong \clh \oplus \clh , $ and ${\dlt}_0= d_1 
\oplus
 d_2 , $ then ${\cll }_0 = -\frac {1} {2} ({d_1}^2 + {d_2}^2 ), \cld 
=
 {\ath }^{\infty } $ satisfies all the properties (i) - ( iii) and 
one
 can construct a quantum stochastic flow driven by $( \pi , {\dlt}_0 
,
 {\cll }_0 ). $ In analogy, one can have the perturbed triple $( \pi 
,
 \dlt , {\cll}) $ where $ \dlt= {\dlt }_1 \oplus {\dlt}_2 $ with $
 {\dlt }_1 = d_1 + d_{r_1} $ and $ {\dlt}_2 = d_2 + d_{r_2} $ and
 $\cll = -\frac {1} {2} ( {\dlt_1 }^2 + {\dlt_2 }^2), \cld =
 {\ath}^{\infty }.$

  Thus we have two triples $(\pi , {\dlt}_0 , 
{\cll}_0)$ and $(\pi , {\dlt} , {\cll})$ both satisfying (i)-(iii). 
 Hence they should give rise to two quantum stochastic processes 
and that they indeed do so is the content of theorem 2.1. Therefore from 
the quantum stochastic point of view also, the two "Laplacians" 
${\cll}_0$ and $\cll$ are equally good candidates for driving the 
processes. Then the question arises: can we associate the same 
geometric features with these two Laplacians or are there 
geometrically discernible changes as we go from the Laplacian 
${\cll_0}$  to the perturbed one $\cll$ ? This will be addressed in 
the following section. \bthm (i)
 The quantum stochastic differential equation (q.s.d.e) $ \cite 
{KRP}$
 for $ x \in {\ath }^{\infty } $ \bea
\label {5}
d j_t^0 (x) = j_t^0 ( i d_1 (x) ) dw_1(t) + j_t^0 ( i d_2 (x) )
dw_2(t) +
j_t^0 ( {\cll }_0 (x))  dt ;{\nonumber}\\
j_0^0 (x) = x \otimes I \eea
has unique solution $j_t^0 $ which is a $*$-homomorphism from $\ath $ 
to \\
$\ath \otimes \clb ( \Gamma ( L^2 ( R_+ ) \otimes C^2 )).$ In fact
$j_t^0 (x) = \al_{( exp2 \pi i w_1 (t),exp2 \pi i w_2 (t))} (x),$
where $( w_1,w_2 ) (t ) $ is the standard two dimensional Brownian
motion. Also \\ $ E j_t^0 (x) = e^{t \cll_0 } (x) $, where $E$ is the
vacuum expectation in the
Fock space $ \Gamma  ( L^2 ( R_+ ) \otimes  C^2  ) $.  \\
(ii) The q.s.d.e in $ \clh \otimes \Gamma $: \bea
\label {6}
d U_t = \sum_{l=1}^2 U_t\{ i {j_t}^0 (r_l ) d {A_l}^\dagger +
i {j_t}^0 ({r_l}^*  )  d {A_l}  - \frac {1} {2} {j_t}^0 ( {r_l}^* r_l  
)  dt  \}  ,{\nonumber}\\
U_0=I \eea has a unique unitary solution $ \cite {EH}.$ Setting $ 
{j_t}
(x) = U_t j_t^0 (x) U_t^*, $ one has the q.s.d.e : \bea
\label{700}
d j_t (x) = \sum_{l=1 }^2 \{ j_t ( i {\dlt }_l (x) ) d {A_l}^{\dagger
  } + j_t ( i {{\dlt }_l}^{\dagger} (x) ) d {A_l} \} + j_t( \cll (x))
dt, \eea and $ E j_t (x) = e^{t \cll } (x).$ \ethm We do not give the
proof here since most of it is available in the references cited
above.  \newsection {Weyl Asymptotics for $\ath$ } For classical
compact Riemannian manifold (M,g) of dimension d with metric g, one
has the natural heat semigroup $\clt_t$ as the expectation semigroup
of the Brownian motion on the manifold \cite{Ros} so that the
Laplace-Beltrami operator $\Delta $ is the generator of $\clt_t$. It
is known \cite{Ros} that $\clt_t $ is an integral operator on $L^2 (
M,dvol ) $ with a smooth integral kernel $\clt_t ( x,y) $, which
admits an asymptotic expansion as $ t \rightarrow 0+ $: \bea
\label{7}
\clt_t(x,y)= \sum_{j=0}^\infty \clt^{(j)} ( x,y) t^{-d/2 +j}, \eea
and that$$vol(M)= \int_M \clt^0(x,x) dvol (x) $$
$$= \lim_{t \rightarrow 0+} t^{d/2} \int_M \clt_t (x,x) dvol (x) =
\lim_{t \rightarrow 0+ } t^{d/2} (Tr \clt_t ), $$
where we have taken
the trace in $ L^2 ( M, dvol )$. Similarly the scalar curvature $s$ 
at
$x \in M$ is given as $s(x)=\frac {1} {6} \clt^{(1)} (x,x) $. This
gives the integrated scalar curvature $$s= \int_M s(x) dvol(x) =\frac
{1} {6} \int_M \clt^{(1)} (x,x) dvol (x)$$
$$
= \frac {1} {6} \lim_{t \rightarrow 0+ } t^{d/2-1} \int [
\clt_t(x,x) - t^{-d/2} \clt^0 ( x,x)] dvol (x)$$
$$
= \frac {1}{6} \lim_{t \rightarrow 0+ } t^{d/2-1} [ Tr \clt_t -
t^{-d/2} vol (M) ]$$
For the non-commutative d-torus ( with d even )
one possibility is to define its volume $V$ and integrated scalar
curvature $s$ by analogy from their classical counterparts as :  \bea
\label {8v}
V(\ath ) \equiv V \equiv \lim_{t \rightarrow 0+ } t^{d/2}Tr \clt_t ,
\eea \bea
\label {8s}
s(\ath ) \equiv s \equiv \frac {1} {6} \lim_{t \rightarrow 0+ }
t^{d/2-1} [ Tr \clt_t -t^{-d/2} V ] \eea where the heat semigroup
$\clt_t$ in the classical case is replaced by the expectation
semigroups of the last section: $\clt_t^0 = e^{t {\cll}_0} $ and the
perturbed one $\clt_t=e^{t \cll} $ respectvely acting on $ L^2 ( \ath
, \tau ) $.  Before we can compute these numbers, we need to study 
the
operators $\cll_0 $ and $ \cll $ in $L^2 ( \tau ) $ more carefully .
The next theorem summarizes their properties for $d=2$ and we have
denoted by $\clb_p$ the Schatten ideals in $\clb ( \clh )$ with the
respective norms.  \bthm
\label {zz}
(i) $ {\cll}_0 $ is a negative selfadjoint operator in $L^2 ( \tau ) 
$
with compact resolvent. In fact $ \cll_0 ( U^m V^n ) = - \frac {1}{2}
( m^2 + n^2 ) U^m V^n ; m,n \in Z $ so that $ {(\cll_0 - z) }^{-1} 
\in
{\clb }_p ( L^2 ( \tau ) ) $ for $ p > 1 $ and $ z \in \rho
(\cll_0)$ \\
\\
(ii) If $ r_1 , r_2 \in \ath^\infty $ and are selfadjoint,  then
$ \cll = \cll_0 + B + A,
$ where $ B = - \frac {1} { 2} ( d_{r_1}^2 + d_{ r_2} ^2 + d_{d_1 (
  r_1)} + d_{d_2 ( r_2) } )$ and $ A = - d_{r_1} d_1 -d_{ r_2} d_2 $,
so that A is compact relative to $\cll_0$ and $\cll $ is selfadjoint
on $ \cld ( \cll_0 ) $ with compact resolvent.

If $r_1,r_2 \in \ath ,$ then $ -\cll = -\cll_0 -B -A $ as quadratic
form on $ D(( -\cll_0)^{\frac {1} {2}}) $ and \bea
\label {9}
{( -\cll + n^2 )}^{-1} = {( - \cll_0 + n^2 )}^{-\frac {1} {2}} {(
  I+Z_n )}^{-1} {( -\cll_0 + n^2 )}^{-\frac {1} {2}} \eea
where \\
$ Z_n = {( - \cll_0 + n^2 )}^{-\frac {1} {2}} (B+A) {( - \cll_0 + n^2
  )}^{-\frac {1} {2}} $, is compact for each n with $B= - \frac {1}
{2} (d_{r_1}^2 + d_{r_2}^2 ), A = \frac {1}{ 2} (d_1 d_{r_1} +d_{r_1}
d_1 +d_2 d_{r_2} +d_{r_2} d_2). $ This defines $\cll $ as a
selfadjoint operator in $L^2 ( \tau ) $ with compact resolvent.
Furthermore, in both cases of (ii) , the difference of resolvents $ 
{(
  \cll - z )}^{-1}-{( \cll_0 - z )}^{-1} $ is trace class for $ z \in
\rho ( \cll ) \cap \rho ( \cll_0 ) $.  \ethm
{\it Proof:--}\\
The proof of (i) is obvious and hence is omitted.  (ii) It is easy to
verify that $\cll = \cll_0 +B +A $ on $\ath^\infty $ and that $A {(
  -\cll_0 + n^2 )}^{-1} $ is compact for every $ n=1,2, \ldots $.
Therefore $( \cll - \cll_0 ) {( - \cll_0 +n^2 )}^{-1} = ( \cll -
\cll_0 ) {( - \cll_0 +1 )}^{-1} ( \cll_0 +1 ) {( - \cll_0 +n^2 
)}^{-1}
\rightarrow 0 $ in operator norm as $ n \rightarrow \infty $. By the
Kato-Rellich theorem \cite{RS}, $\cll $ is selfadjoint and since \\
${(- \cll +n^2 )}^{-1} = {(- \cll_0 + n^2)}^{-1} {[ 1 + ( \cll_0 -
  \cll ) {(- \cll_0 + n^2 )}^{-1} ]}^{-1} $ for sufficiently large
$n$, one also concludes that $\cll $ has compact resolvent.
Furthermore for $ z \in \rho ( \cll ) \cap \rho ( \cll_0) ,$ $$
{(
  \cll -z)}^{-1} - {( \cll_0 - z )}^{-1} = {( \cll_0 - z )}^{-1} [ 1 
+
(\cll - \cll_0 ) {( \cll_0 - z )}^{-1} ]^{-1} ( \cll_0 - \cll ){(
  \cll_0 - z )}^{-1} $$
Since $ ( \cll - \cll_0 ){( -\cll_0 +
  n^2)}^{-\frac {1} {2}} $ is bounded , ${( -\cll_0 + n^2)}^{-\frac
  {1} {2} } \in \clb_3 ( L^2 ( \tau )) $ and since ${( -\cll_0 + z
  )}^{-1} \in \clb_{3/2} ( L^2 ( \tau )),$ It follows that $ {( \cll 
-
  n^2)}^{-1}-{( \cll_0 - n^2 )}^{-1} $ is trace class for $n=1,2,
\ldots $ by the Holder inequality.

When $ r_1, r_2 \in \ath ,$ we cannot write the expression for $\cll 
$
as above on $\ath^\infty$, since $ r_1, r_2 $ may not be in the 
domain
of the derivations $ d_1 , d_2 $. For this reason, we need to define 
$
-\cll $ as the sum of quadratic forms and standard results as in 
\cite
{RS} can be applied here . From the structure of $B$ and $A$ it is
clear that $Z_n$ is compact for each $n$ and hence an identical
reasoning as above would yield that $ \| Z_n \| \raro 0 $ as $ n 
\raro
\infty $ and therefore $ {(I+ Z_n)}^{-1} \in \clb $ for sufficiently
large $n$ and the right hand side of (\ref {9}) defines the operator
$- \cll $ associated with the quadratic form with $ D ((- \cll
)^{\frac {1} {2}} ) = D ((- \cll_0 )^{\frac {1} {2}} ). $ Clearly $$
{( -\cll + n^2 )}^{-1} - (-\cll_0 +n^2 )^{-1}= -{( - \cll_0 + n^2
  )}^{-\frac {1} {2}} {( I+Z_n )}^{-1} Z_n {( -\cll_0 + n^2 
)}^{-\frac
  {1} {2}} $$
$$
= - (- \cll_0 +n^2 )^{-\frac {1} {2}} (I+Z_n)^{-1} (-\cll_0 +n^2
)^{-\frac {1} {2}} (B+A) ( - \cll_0 + n^2 )^{-1}$$
for sufficiently
large $n$ and since
$$
(- \cll_0 + n^2 )^{-\frac {1} {2}} \in \clb_3, (- \cll_0 +n^2
)^{-\frac {1} {2}} A (- \cll_0 +n^2 )^{-\frac {1} {2}} \in \clb_3,$$
it is clear that $ {( \cll - n^2)}^{-1}-{( \cll_0 - n^2 )}^{-1} $
is trace class.   \qed     \\
The next theorem studies the effect of the perturbation from $\cll_0 
$
to $\cll$ on the volume and the integrated sectional curvature for
$\ath$.  \bthm (i) The volume V of $\ath (d=2) $ as defined in (\ref
{8v}) is invariant under the perturbation
from $\cll_0$ to $\cll$.  \\
(ii) The integrated scalar curvature for $r \in \ath^\infty$,  in general is not invariant under the above perturbation.

\ethm {\it Proof :-- } We need to compute $Tr ( e^{t \cll } - e^{t
  \cll_0 } ).$ Note that if $ r_1, r_2 \in \ath^\infty,$ then $e^{t
  \cll } - e^{t \cll_0 } = - \int_0^t e^{(t-s) \cll } ( \cll -\cll_0 
)
e^{s \cll_0 } ds $ which on two iterations yields:
$$e^{t \cll } - e^{t \cll_0 } = - \int_0^t e^{(t-s) \cll_0 } ( \cll
-\cll_0 ) e^{s \cll_0 } ds +\int_0^t dt_1 e^{(t-t_1 ) \cll_0 } (\cll
-\cll_0 )\times$$
$$\int_0^{t_1} dt_2 e^{(t_1 -t_2 ) \cll_0 } ( \cll - \cll_0 ) e^{ t_2
  \cll_0 } - \int_0^t dt_1 e^{( t- t_1) \cll } ( \cll -\cll_0 )
\int_0^{t_1} dt_2 e^{(t_1 -t_2 ) \cll_0 } (\cll- \cll_0 )\times $$
\bea
\label{10}
\int_0^{t_2} dt_3 e^{( t_2 -t_3 ) \cll_0 } ( \cll - \cll_0 ) e^{ t_3
  \cll_0 } \equiv I_1 (t) + I_2 (t) + I_3 (t).  \eea For estimating
the trace norms of these terms , we note that the $\clb_p $-norm of
$(\cll - \cll_0 ) e^{s \cll_0 } $ is estimated as $$
\| ( \cll -\cll_0
) e^{s \cll_0 } \|_p = \| ( B+ A ) e^{s \cll_0} \|_p \leq \| B \| \|
e^{s \cll_0 } \|_p + $$
$$
c_1 (\| d_1 e^{s \cll_0} \|_p + \| d_2 e^{s \cll_0 } \|_p ) \leq
c^{\prime \prime } ( \| e^{s \cll_0 } \|_p + \| d_2 e^{s \cll_0 } 
\|_p
)$$
$$
\leq c' (s^{-{p^{-1}}} + s^{-{p^{-1}}-\frac {1} {2}} ) \leq c ~
s^{-{p^{-1}}-\frac {1} {2}} $$
for constants $c,c_1,c^\prime
,c^{\prime \prime} $ since we are interested only for the region \\ $
0 < s \leq t \leq 1 .$ Using Holder inequality for Schatten norms and
the fact that $$\| ( \cll - n^2 )^{-1} \| \leq \| ( \cll_0 - n^2
)^{-1} [ 1 + (\cll -\cll_0) (\cll_0 - n^2 )^{-1} ]^{-1} \| \leq \frac
{2} { n^2 } $$
for sufficiently large $n$.  We get for the third term
in \ref {10}
$$
\| I_3 (t) \|_1 \leq 2 \int_0^t dt_1 \int_0^{t_1} dt_2 \| ( \cll -
\cll_0 ) e^ {( t_1 - t_2 ) \cll_0 } \|_{p_1} \times $$
$$
\int_0^{t_2} dt_3 \| (\cll - \cll_0 ) e^{(t_2 - t_3 ) \cll_0 }
\|_{p_2} \| ( \cll -\cll_0 ) e^{t_3 \cll_0 } \|_{p_3}
$$
$$
\leq c(p_1,p_2,p_3) \int_0^t t_1^{-\frac {1} {2}} dt_1 \raro 0
$$
as $ t \raro 0$ where $p_1^{-1} + p_2^{-1} + p_3^{-1} =1 $. A very
similar estimate shows that
$$
\| I_1(t) \|_1 \leq \int_0^t ds \| e^{(t-s) \cll_0} \|_{p_1} \|
(\cll- \cll_0 ) e^{s \cll_0 } \|_{p_2} \leq c t^{-\frac {1} {2}}$$(
with $ p_2 > 2 ~ and ~ p_1^{-1} + p_2^{-1}=1)$ and $$
\| I_2 (t) \|_1
\leq \int_0^t dt_1 \| e^{( t - t_1 ) \cll_0 } \|_{p_1} \int_0^{t_1}
dt_2 \| (\cll - \cll_0 ) e^{(t_1 - t_2 ) \cll_0 } \|_{p_2} \| ( \cll 
-
\cll_0 ) e^{t_2 \cll_0 } \|_{p_3}\leq c^{\prime}, $$
(with $ p_1^{-1}
+ p_2^{-1} + p_3^{-1} =1, $ in particular the choice $ p_1=p_2=p_3=3 
$
will do) a constant independent of t . From this it follows that \\ $
\lim_{t \raro 0+ } t~ Tr (e^{t \cll} - e^{t \cll_0} ) =0 $ and thus
the invariance of volume under perturbation.

In the case when $r_1,r_2 \in \ath $ only , then $ \cll - \cll_0 = B 
+
d_1 B_1 + d_2 B_2 + B_1^{\prime} d_1 + B_2^{\prime } d_2$ where $
B,B_1,B_1^{\prime}, B_2 , B_2^{\prime} $ are bounded. Therefore the
term like $e^{(t-s) \cll_0 } d_1 B_1 e^{s \cll_0} =[ e^{s \cll_0}
B_1^* d_1 e^{(t-s) \cll_0 } ]^*$ admits similar estimates
as above and the same result  follows. \\
(ii) From the expression (\ref {8s}) for the integrated scalar
curvature $s$, we see
that for $d=2$  \\
\bea
\label {20}
s (\cll) - s(\cll_0) = \frac {1} {6} ~ {\lim_{t \raro 0+}} Tr (e^{t
  \cll}- e^{t \cll_0} ) \eea if it exists, and conclude that the 
contribution
   to (\ref {20}) from the term $I_3(t)$ vanishes as we have seen in 
(i).We
claim that though ${\| I_2(t) \|}_1 \leq$ constant, $TrI_2(t) \raro 
0$
as $ t \raro 0+ $. In fact since the integrals in $I_2(t) $ converges
in trace norm
$$Tr I_2(t) = \int_0^t dt_1 \int_0^{t_1} dt_2 Tr ( (\cll - \cll_0
)e^{(t_1-t_2) \cll_0}
(\cll -\cll_0)  e^{(t - t_1 +t_2 ) \cll_0 } )$$
and by a change of variable we have that \\
$| Tr I_2(t) | \leq t \int_0^t\|(\cll -\cll_0) e^{s \cll_0} (\cll
-\cll_0) e^{(t-s) \cll_0} \|_1 ds $ For $ r \in \ath^\infty $, the
perturbation $(\cll -\cll_0)$ is of the form $ b_0 + b_1 d_1 + b_2 
d_2
$ with $ b_i \in \clb ( \clh ) $ for $ i=0,1,2$ and the
Hilbert-Schmidt norm estimates are as follows : $$
\| (\cll -\cll_0)
e^{s \cll_0} \|_2 \leq \| b_0 \| \| e^{s \cll_0} \|_2 + \sqrt{2} (
\|b_1 \| + \| b_2 \| ) \| {(- \cll_0)}^{\frac {1} {2}} e^{s \cll_0}
\|_2 \leq c( s^{-\frac {1} {2}} + s^{-\frac {3}{4}}).  $$
Therefore
$$| Tr I_2(t) | \leq c t \int_0^t ( s^{-\frac {1} {2}} + s^{-\frac
  {3}{4}} )  ( (t-s)^{-\frac {1} {2}}
+(t-s)^{-\frac {3}{4}} ) $$
and this clearly converges to zero as  $ t \raro 0+ $. \\
This leaves only $I_1(t)$ contribution so that $$
6(s(\cll) -
s(\cll_0)) = - \lim_{t \raro 0+ } t~ Tr ( ( \cll -\cll_0) e^{t 
\cll_0}
) .$$
As before we note that $ ( \cll - \cll_0 ) $ contains two kinds
of terms :\\ $ B= -\frac {1} {2} ( d_{r_1}^2
+  d_{r_2}^2 ) , A = -\frac {1} {2} ( d_{r_1}  d_1 + d_1 d_{r_1} + 
d_{r_2} d_2 + d_2 d_{r_2} )$  \\
We show that the term $Tr (A e^{t \cll_0} )=0 $ for all $t >0$. It
suffices to show that $ Tr ( d_r d_1 e^{t \cll_0}) =0 $ for $ r \in
\ath^\infty $ and for this we note that
$$
Tr ( d_r d_1 e^{t \cll_0}) = \sum_{m,n}  < U^m V^n , d_r d_1
e^{t \cll_0} ( U^m V^n ) > $$
$$
= \sum_{m,n} m e^{-t/2 (m^2 + n^2 )} \tau ( V^{-n} U^{ -m} d_r (
U^m V^n )) = \sum_{m,n} m e^{-t/2 (m^2 + n^2 )} \tau ( V^{-n} U^{-m} 
r
U^m V^n -r ) =0 $$
identically.  This leaves only the contribution due
to $B$. Thus \bea
\label {sl}
s(\cll ) - s(\cll_0)= \frac {1}{12} \lim_{t \raro 0+} t ~ Tr
((d_{r_1}^2 + d_{r_2}^2 )e^{t \cll_0} ), \eea if it exists.
However since $\{t Tr (({d_{r_1}}^2 +{d_{r_2}}^2) e^{t {\cll}_0}) \} $ is bounded as $ t \raro 0+$, we can and will interpret the above limit as a special kind of Banach limit  as in  Connes [2,p.563] \bea
 s (\cll)- s({\cll}_0 ) = \frac {1}{12} {Lim}_{t^{-1} \raro \omega} tTr (({d_{r_1}}^2 +{d_{r_2}}^2) e^{t {\cll}_0}) \\ = \frac {1}{12} Tr_{\omega}
 (({d_{r_1}}^2 +{d_{r_2}}^2 ) {\hat {\cll_0}}^{-1}  \eea
The notation ${\hat {\cll_0}}$ will be explained in the next section.
 In the following
 we show that in general the right hand side of (3.8) is strictly positive.

 For example  set $r_1 = ( U+ U^{-1}) $ and $r_2=0$, then
 $ r_1,r_2
\in \ath^\infty,$ and $$6( s(\cll) -s(\cll_0) )= \frac {1} {2} 
{Lim}_{t^{-1} \raro \omega} t \sum_{m,n} e^{-t/2 (m^2 + n^2 )}  < U^m V^n ,
d_{r_1}^2 ( U^m V^n ) > $$
$$
= 2^{-1} {Lim}_{t^{-1} \raro \omega} t \sum_{m,n} e^{-t/2 (m^2 + n^2 )} \tau
( (1- \lmd^{-n} )^2 \lmd^{2n} U^2 + (1 - \lmd^n)^2 \lmd^{-2n} U^{-2} 
+
(2- \lmd^n -\lmd^{-n} ))$$
$$
=2^{-1} {Lim}_{t^{-1} \raro \omega} t \left( 2 \sum_{m=1}^{\infty} e^{-m^2 t 
/2}
  +1 \right ) \left( 8 \sum_{n=1}^\infty {sin}^2 (\pi \theta n)
  e^{-n^2 t/2} \right) $$
Next note that for $ 0 < t < 2$  $$
\sqrt{t} \sum_{n=1}^{
  \infty}{sin}^2 (\pi \theta n) e^{-n^2 t/2} \geq \sqrt{t}
\sum_{n=1}^{ [ \sqrt{2/t}]} {sin}^2 (\pi \theta n) e^{-n^2 t/2} $$
$$\geq e^{-1} (\sqrt {2} - \sqrt{t} ) \sum_{n=1}^{ [ \sqrt{2/t}]} [
\sqrt{2/t}]^{-1} {sin}^2 \pi (n \theta - [n \theta] ) = e^{-1} (\sqrt
{2} - \sqrt{t} ) E ({sin}^2 \pi X_t ),$$
where for each $0 < t \leq
2,~~ X_t $ is a $[0,1]$-valued random variable with Probability$(X_t 
=
k \theta -[k \theta] )= [\sqrt{2/t}]^{-1} $ for $k=1,2,\ldots , 
[\sqrt
{\frac {2} {t}}]$ and $E$ is the associated expectation. Since $
\theta $ is irrational, it is known that (\cite{Hel}) as $t \raro 
0+$,
the random variable $X_t$ converges weakly to one with uniform
distribution on $[0,1]$ and therefore $$\liminf_{t \raro 0+}\sqrt{t}
\sum_{n=1}^{ \infty}{sin}^2 (\pi \theta n) e^{-n^2 t/2} \geq \lim_{t
  \raro 0+}\sqrt{t} \sum_{n=1}^{ [ \sqrt{2/t}]} {sin}^2 (\pi \theta 
n)
e^{-n^2 t/2} $$
$$
\geq \sqrt{2} e^{-1} \int_0^1 {sin}^2 \pi x dx =
{(\sqrt{2}e)}^{-1}.  $$
We also have by Connes (page 563) \cite {Con} $\lim_{t \raro 0+} \sqrt{t}\sum_{m=1}^{\infty} e^{-m^2 t/2} =\frac {\sqrt{\pi}}{\sqrt{2}}$
Now by the general properties of the limiting procedure as expounded in
\cite {Con} $$s (\cll )- s ( \cll_0 ) \geq \frac {2 \sqrt{\pi}}{3e} $$
\qed

{\bf Remark}:--  From the expression for $s( \cll_0) $, we see 
that
for $d=2$,  $s(\cll_0 ) =\lim_{t \raro 0+} ( Tr e^{t \cll_0} - \frac 
{V}
{t}) $. Since the expression for $ Tr e^{t \cll_0} $and the volume 
$V$
are exactly the same as in the case of classical two-torus with its
heat semigroup, the integrated scalar curvature for $\cll_0$ is the
same as in the classical case, which is clearly zero. Therefore $ s(\cll)$ is strictly positive for the case considered here.

\newsection{ Spectral Triple on $\ath^\infty,$ its perturbation and
  cohomology } Following Connes \cite {Con} we consider the even
spectral triple \\ $(\cla= \ath^\infty, \clh = L^2(\tau) \oplus L^2 (
\tau ) , D_0 , \Gamma ) $ where $D_0$ , the unperturbed Dirac
operator=
$  \left( \begin{array} {cc}  0 & d_1+id_2  \\
    d_1 - id_2 & 0
\end {array} \right) \equiv i \gma_1 d_1 (a) +i \gma_2 d_2 (a)$ in 
$\clh$.
Here $ \gma_1,\gma_2$ are the $2 \times 2 $ clifford matrices. The 
grading
operator is given by $ \Gamma =  \left( \begin{array}
{cc}  I & 0  \\
                               0   & -I
\end {array} \right)  .$ One easily verifies that $ a \Gamma = \Gamma 
a , \Gamma^*= \Gamma =
\Gamma^{-1}, \Gamma D_0 = -D_0  \Gamma$. Note also $D_0$ has compact 
resolvent since
$D_0^2 = -2  \left( \begin{array} {cc}  \cll_0& 0  \\
                               0   & \cll_0
\end {array} \right)  $  and $ ker D_0 = ker \cll_0 \otimes C^2$ is 
two dimensional.
The perturbed spectral triple is taken to be $(\cla, \clh,D, \Gamma)$ 
where $ D= D_0 +
 \left( \begin{array} {cc}  0 & d_r  \\
                               d_{r^*}   & 0
\end {array} \right) $ for some $r \in \ath^\infty$.  It is not 
difficult to see that
$D_0 $ and $D$ are both essentially selfadjoint on $ \cla  \subseteq 
\L^2(\tau) $ and that the
 perturbed triple is also an even one. Here, as in Connes \cite 
{Con}, by the volume form $v(a) $ on
$\cla $ we mean the linear functional $v(a)= \frac{1} {2} {Tr}_w (a | 
\hat {D}  |^{-2}  P) $
where ${Tr}_w  $ is the Dixmier trace \cite {Con}, and we have used 
the notation that for a
selfadjoint operator $T$ with compact resolvent $\hat {T} = T 
|_{{N(T)}^\perp} \equiv TP$, where $P$ is the
projection on ${N(T)}^{\perp}$.
Next we prove that the volume form is invariant under the above 
perturbation.
 For this we need a lemma.
\blmma
Let T be a selfadjoint operator with compact resolvent such that 
${\hat {T}}^{-1} $
is Dixmier trace-able. Then for $a \in \cla $ and every $z \in \rho 
(T),$ \\ $ {Tr}_w
( a {\hat {T}}^{-1} P) = {Tr}_w ( a {(T-z)}^{-1}).$
\elmma
{\it Proof:-- } \\
Note that $ {(T-z)}^{-1} = {(\hat {T} - z )}^{-1} P \oplus -z^{-1} 
P^{\perp}$ and  $ P^{\perp}$
is finite dimensional. Therefore $ {Tr}_w ( a {(T-z)}^{-1}) = {Tr}_w 
(PaP {( \hat{T} -
z )}^{-1} P).$ On  the other hand\\ $ {Tr}_w ( PaP {\hat{T}}^{-1} P - 
PaP {(\hat{T} -z )}^{-1} P )
=-z {Tr}_w ( PaP {\hat {T}}^{-1} {(\hat {T} -z )}^{-1} P ) =0 $, 
since ${\hat {T}}^{-1} $ is
Dixmier trace-able and ${(\hat{T} -z )}^{-1}  $ is compact \cite 
{Con}\qed
\bthm
If we set $v_0 (a) = \frac {1} {2} {Tr}_w ( a | {\hat {D}}_0|^{-2} ) 
$ and
$v(a)= \frac {1} {2} {Tr}_w (a {| \hat {D} |}^{-2}) $ for $a \in \cla 
$, then $ v_0(a)=v(a)$
\ethm
{\it Proof :--} \\
Note that $ D^{2}=-2   \left( \begin{array} {cc}  \cll_1 & 0  \\
                              0  & \cll_2
\end {array} \right) ,$ where \\ $ \cll_1= \cll_0 + d_r d_{r^*}  + ( 
d_1 d_{r^*} + d_r d_1 )
+ i ( d_2 d_{r^*} - d_r d_2 )$ and \\ $ \cll_2 = \cll_0 + d_{r^*} d_r  
+ ( d_1  d_r  + d_{r^*}
d_1 ) +i ( d_2 d_{r^*} -d_r d_2 ) ,$ and that by theorem (\ref {zz}) 
of section 3, both $\cll_1$
and $\cll_2$ have compact resolvents with $  P_1, P_2 $ projections 
on ${\cln (\cll_1)}^\perp$
and  ${\cln (\cll_2)}^\perp$ respectively.
Therefore by the previous lemma for $ Im z \neq 0$ $$ v(a) = {Tr}_w ( 
a {(-\hat {\cll_1})}^{-1} P_1 ) +
{Tr}_w ( a {(-\hat {\cll_2})}^{-1} P_2 )$$  $$ = {Tr}_w ( a{(- \cll_1 
-z )}^{-1} +
a{(- \cll_2 -z )}^{-1}) $$  $$= {Tr}_w ( a{( -\cll_0 -z )}^{-1}
+a{( -\cll_0 -z )}^{-1}) +{Tr}_w ( a{( -\cll_1 -z )}^{-1}$$  $$
-a{( -\cll_0 -z )}^{-1})+{Tr}_w ( a{( -\cll_2 -z )}^{-1}
-a{( -\cll_0 -z )}^{-1}) = v_0(a) $$
since $ {( -\cll_i -z )}^{-1} -  {( -\cll_0 -z )}^{-1} $ is trace 
class for $ i=1,2 $
\qed
\vspace {.5 in}

We say that two spectral triples $(\cla_1, \clh_1, D_1)$ and 
$(\cla_2, \clh_2, D_2)$ are unitarily equivalent
 if there is a unitary operator $U : \clh_1 \raro \clh_2$ such that 
$D_2=U D_1 U^*$ and $\pi_2(.)=U \pi_1(.) U^*,$
 where $\pi_j,~j=1,2$ are the representation of $\cla_j$ in $\clh_j$ 
respectively. Now, we want to prove that in general
 the perturbed spectral triple is not unitarily equivalent to the
  unperturbed one. Let $\Omega^1(\ath^\infty)$ be the
 universal space of $1$-forms (\cite{Con}) and $\pi$ be the 
representation of
  $\Omega^1 \equiv \Omega^1(\ath^\infty)$
 in $\clh$ given by
$$ \pi(a)=a, \pi(\dlt(a))=[D,a],$$
where $ \dlt $ is the universal derivation. \\
Note that  $[D,a]=i [ \dlt_1(a)\gma_1+\dlt_2(a)\gma_2 ],$ where  
$r_1={\rm
Re}~r,~r_2={\rm Im}~r ~\dlt_1= d_1+d_{r_1}, \dlt_2= d_2 + d_{r_2}.$

\bthm
(i) Let $r=U^m$, then
 $\Omega^1_D(\ath^\infty) := \pi(\Omega^1)=\ath^\infty \oplus 
\ath^\infty.$\\
(iii) $\Omega^2(\ath^\infty)=0$ for $r=U^m$.
\ethm
{\it Proof :-}\\
(i) Clearly $\pi(\Omega^1) \subseteq \ath^\infty \gma_1+\ath^\infty 
\gma_2.$
The other inclusion follows from the facts that
 $\dlt_2(U^k)=0, \dlt_1(U^k)$ is invertible, and that
 $\dlt_2(V^l)$ is invertible for sufficiently large $l$.\\

(iii) Let $J_1={\rm Ker} \pi|_{\Omega^1},~J_2={\rm Ker} 
\pi|_{\Omega^2}.$ Then $J_2+\dlt J_1$ is an ideal, implying
 that $\pi(\dlt J_1)=\pi(J_2+\dlt J_1)$ is a nonzero submodule of 
$\pi(\Omega^2)\subseteq \ath^\infty \oplus \ath^\infty$.
 Since $\ath^\infty$ is simple there are two possibilities, namely 
either $\pi(\dlt J_1) \cong \ath^\infty$, or
 $\pi(\dlt J_1) =\ath^\infty \oplus \ath^\infty.$
  To rule out the first possibility we take a closer look at $ J_1 $ 
and
  $\pi ( \dlt J_1)$. $J_1 = \{ \sum_i a_i \dlt (b_i) |  \sum_i a_i 
\dlt_1 (b_i)=0,
   \sum_i a_i \dlt_2 (b_i) =0 \} $. Using the fact that $\dlt_1, 
\dlt_2$
    are derivations we get
    \bea \label{200} \sum_i \dlt_1(a_i) \dlt_2 (b_i) =-
    \sum_i a_i \dlt_1 (\dlt_2 (b_i))  \\
     \label {201} \sum_i \dlt_2(a_i) \dlt_1 (b_i) =- \sum_i a_i 
\dlt_2(\dlt_1 (b_i))
     \eea
      for $  \sum_i a_i \dlt (b_i) \in J_1 $
     $$\pi ( \sum_i \dlt (a_i) \dlt ( b_i) )=
      \sum_i (\dlt_1(a_i) \gma_1 + \dlt_2 (a_i) \gma_2 ) ( 
\dlt_1(b_i) \gma_1
      +\dlt_2 (b_i) \gma_2 ) $$
  $$ =\sum_i ( \dlt_1 (a_i) \dlt_1(b_i) + \dlt_2 (a_i) \dlt_2 (b_i))
      + \sum ( \dlt_1 (a_i) \dlt_2(b_i)- \dlt_2 (a_i) \dlt_1(b_i)) 
\gma_{12},$$
     where  $ \gma_{12}=\gma_1 \gma_2=- \gma_2 \gma_1.$
      Taking $x = U^{-1} \dlt (U) + U \dlt (U^{-1}) \in \Omega^1 $
      it is easy to verify that $ x \in J_1 $ and $ \pi (\dlt x ) = 
-2 $.
     This proves $ \ath^\infty \oplus 0 \subseteq \pi ( \dlt J_1) $. 
We show that
     the inclusion is proper by showing  the nontriviality of
     coefficient of $ \gma_{12}$.
      Using \ref {200}, \ref {201} we get coefficient of $ \gma_{12}$ 
to be
    $  \sum a_i [\dlt_1,\dlt_2] (b_i)= \sum -im a_i [r_1 , b_i].$
       As before we can find   $n_0$  such that for $l \geq n_0, ~
       \dlt_2 (V^l) $ is invertible. If we now choose $ a_1 = I , b_1 
= V^{n_0},
       a_2 = - \dlt_2 (V^{n_0}) \dlt_2 (V^l)^{-1} , b_2= V^l,
       a_3=(-a_1 \dlt_1(b_1) - a_2 \dlt_2(b_2) ) U^{-1}, b_3=U $,
       then the vanishing of the coefficient of $\gma_{12}$ will
imply that $[ r_1 , V^{n_0} ]= \dlt_2 (V^{n_0} ) \dlt_2 (V^l)^{-1}
 [ r_1, V^l]$ for all $ l \geq n_0$ and we note that while the left 
hand side is nonzero and independent of l, the right hand side 
converges to 0 as $l \raro \infty $ leading to a contradiction. 
Therefore
  $\ath^\infty \oplus \ath^\infty = \pi (\dlt J_1) \subseteq 
\pi(\Omega^2) \subseteq
 \ath^\infty \oplus \ath^\infty.$ Hence 
$\Omega^2_D(\ath^\infty)=\frac{\pi(\Omega^2)}{\pi(\dlt J_1)} =0$.
\qed

Thus we have the following :
\bthm
The spectral triples $(\ath^\infty, \clh, D_0)$ and $(\ath^\infty, 
\clh, D)$ are not unitarily equivalent for $r=U^m$.
\ethm
The proof is clear since $\Omega^2_{D_0}(\ath^\infty)=\ath^\infty \ne 
0 = \Omega^2_D(\ath^\infty).$

Classically there is a correspondence between connection form and 
covariant differentiation. This correspondence
comes from the duality between the module of derivations and the 
module of sections in the cotangent bundle.
Unfortunately there is no such duality in the non-commutative 
context. Here for defining the connection form
we visualize it more as the connection form arising from covariant 
differentiation. We need to do so because
if we take the existing definition~\cite {Fro} then the curvature 
form becomes trivial.

Let $\clk$ be the vector space of all derivations $d :\ath^\infty 
\raro \ath^\infty$. This space is same as
$\{ c_1d_1+ c_2 d_2 +[r,.] : r \in \ath^\infty \}$ for almost all
$ \theta $ ( lebesgue)  \cite {BEJ} for the rest of this section we 
will be using those $\theta 's $ only.
 Let $\dlt_{mn}$ be the element of $\clk$ given by
 $\dlt_{mn}(a)=[U^m V^n, a].$ We turn $\clk$ into an inner product 
space by requiring that $\{ d_1, d_2, \dlt_{mn}\}$
 to be orthonormal, for example as in \cite {JOR}. Let $\cle$ be any 
normed $\ath^\infty$-module. For $\dlt \in \clk$, let $c_\dlt : \cle 
\ot \clk
 \raro \cle,$ be the contraction with respect to $\dlt$. Topologize 
$\cle \ot \clk$ with the weak topology inherited
 from $c_\dlt, \dlt \in \clk$. Then a {\it connection} is a 
complex-linear map $\Nab : \cle \raro \cle \ot \clk$ such
 that $c_\dlt \Nab (\xi a)=c_\dlt \Nab(\xi)a+\xi \dlt(a),~\forall 
\dlt \in \clk.$
\bthm
Suppose that $\Nab_1, \Nab_2$ are maps from $\cle$ to $\cle$ 
satisfying
$$ \Nab_i(\xi a)=\Nab_i(\xi)a+\xi d_i(a),~i=1,2.$$
 Then the map $\Nab$ given  by
$$ \Nab(\xi)=\Nab_1 \ot d_1+\Nab_2 \ot d_2 -\sum \xi U^mV^n \ot 
\dlt_{mn}$$
 is  well-defined and is a connection.
\ethm
{\it Proof :-}\\
Let $\dlt \in \clk$, such that $\dlt=c_1 d_1 +c_2 d_2 +\sum c_{mn} 
\dlt_{mn},$ where $\{c_{mn}\} \in \cls(\IZ^2) \subseteq l_1(\IZ^2)$. 
Therefore the sum in the right hand side of the definition of $\Nab$ 
converges in the topology referred above.
 The rest is straightforward.
\qed

It is clear from the definition of $\Nab $ in the above theorem
that $ \Nab_j = c_{d_j} \Nab $ for $ (j=1,2) $. We also set
$\Nab_r = c_{d_r} \Nab $ for $ r \in \ath^\infty $
\bdfn
Let $R :\clk \ot \clk \raro \cll(\cle)$ be the map given by 
$R(\dlt_1, \dlt_2)=c_{[\dlt_1,
\dlt_2]}  \Nab -[ c_{\dlt_1} \Nab, c_{\dlt_2} \Nab ].$ We call $R$ 
the curvature $2$-form
 associated with the connection $\Nab.$
\edfn
\bthm
We have
$$R(d_1, d_2)=R(d_1+d_{r_1}, d_2+d_{r_2}).$$
\ethm
{\it Proof :-}\\
$ [d_1+d_{r_1},d_2+d_{r_2}]=[d_1(r_2),.]-[d_2(r_1),.]+[ 
[r_1,r_2],.].$
 So we have
 \bean
\lefteqn{R(d_1+d_{r_1},d_2+d_{r_2})(\xi) }\\
&=& -\xi d_1(r_2)+\xi d_2(r_1)-\xi[r_1,r_2]-(\Nab_1+\Nab_{r_1})
(\Nab_2 \xi -\xi r_2)+(\Nab_2+
 \Nab_{r_2})(\Nab_1 \xi-\xi r_1)\\
&=& -[\Nab_1, \Nab_2]\xi +\Nab_1(\xi r_2)+(\Nab_2 \xi)r_1-\xi r_2 
r_1-\Nab_2(\xi r_1) \\
&-& (\Nab_1 \xi) r_2+\xi r_1 r_2 -\xi d_1(r_2)+\xi d_2(r_1)-\xi [r_1, 
r_2]\\
&=& -[\Nab_1, \Nab_2]\xi = R(d_1, d_2)(\xi) ~~({\rm since} ~~[d_1, 
d_2]=0).\\
\eean
\qed

{\bf Remark}:-- In section 3, we have seen that the integrated scalar 
curvature
under the perturbed Lindbladian is different from zero, whereas in 
section 4,the
curvature 2-form has been shown to be invariant under the same 
perturbation.

\newsection{Non-commutative $2d$-dimensional space}
In this section we shall discuss the geometry of the simplest kind of 
noncompact manifolds, namely the
 Euclidean $2d$-dimensional space and its noncommutative counterpart. 
Let $d \geq 1$  be an integer and
 let $\cla_c \equiv C_0(\IR^{2d})$, the (nonunital) $C^*$-algebra of 
all complex-valued continuous functions on $\IR^{2d}$
 which vanish at infinity. Then $\partial_j (j=1,2, \ldots, 2d),$ the 
partial derivative
in the $j$-th direction, can be
 viewed as a densely defined derivation on $\cla_c$, with the domain
 $\cla_c^\infty \equiv C_c^\infty (\IR^{2d}),$
 the set of smooth complex valued functions on $\IR^{2d}$ having 
compact support. We consider the Hilbert space $L^2(\IR^{2d})$
 and naturally imbed $\cla_c^\infty$ in it as a dense subspace. Then 
$i \partial_j$ is a densely defined symmetric linear
 map on $L^2(\IR^{2d})$ with domain $\cla_c^\infty$, and we denote 
its self-adjoint extension by the same symbol.
 Also, let $\clf$ be the Fourier transform on $L^2(\IR^{2d})$ given 
by
$$ \hat{f}(k) \equiv (\clf f)(k)= (2 \pi)^{-d} \int e^{-ik.x} 
f(x)dx,$$
 and $M_\varphi$ be  the operator of multiplication by the function 
$\varphi$. We set $\widetilde{M_\varphi}=\clf^{-1}
 M_\varphi \clf,$ thus $i \partial_j=\widetilde{M_{x_j}}.$ $\Dlt 
\equiv \widetilde{M}_{-\sum x_j^2}$ is the self-adjoint negative
 operator, called the $2d$-dimensional Laplacian. Clearly, the 
restriction of $\Dlt$ on $\cla_c^\infty$ is the differential
 operator $\sum_{j=1}^{2d} \partial_j^2.$ Let $h=L^2(\IR^d)$ and 
$U_\al, V_\beta$ be two strongly continuous groups of
 unitaries in $h$, given by the following :
$$ (U_\al f)(t)=f(t+\al),~~(V_\beta f)(t)=e^{it. \beta} f(t),~~\al, 
\beta, t \in \IR^d,~ f \in C_c^\infty(\IR^d).$$
Here $t.\beta$ is the usual Euclidean inner product of $\IR^d$.
 It is clear that
\bea
\label{weyl}
U_\al U_{\al^\prime}
= U_{\al +\al^\prime},{\nonumber}\\
  V_\beta V_{\beta^\prime}=V_{\beta+\beta^\prime},{\nonumber}\\
 U_\al V_\beta = e^{i \al . \beta} V_\beta U_\al{\nonumber}.\\
\eea
For convenience, we define a unitary operator $W_x$ for $x =(\al, 
\beta) \in \IR^{2d}$ by
$$ W_x=U_\al V_\beta e^{-\frac{i}{2} \al. \beta},$$
 so that the Weyl relation (\ref{weyl}) is now replaced by $W_x 
W_y=W_{x+y} e^{\frac{i}{2} p(x,y)},$ where
 $p(x,y)=x_1.y_2-x_2.y_1,~{\rm for}~x=(x_1,x_2), y=(y_1,y_2).$ This 
is exactly the Segal form of the Weyl relation (\cite{Fol}).
   For $f$ such that $\hat{f} \in L^1(\IR^{2d})$, we set
$$b(f)= \int_{\IR^{2d}} \hat{f}(x)W_x dx \in \clb(h).$$
 Let $\cla^\infty$ be the $\ast$-algebra generated by $\{ b(f) | f 
\in C_c^\infty (\IR^{2d}) \}$ and let $\cla$ be the
 $C^*$-algebra generated by $\cla^\infty$ with the norm inherited 
from $\clb(h)$. It is easy to verify using the commutation
 relation (\ref{weyl}) that
$ b(f)b(g)=b(f \odot g)$ and $b(f)^*=b(f^\natural)$, where
 $$ \widehat{(f \odot g)}(x)=\int 
\hat{f}(x-x^\prime)\hat{g}(x^\prime) e^{\frac{i}{2} p(x, x^\prime)} 
dx^\prime;~~ f^\natural
(x)=\bar{f}(-x).$$ We define a linear functional $\tau$ on 
$\cla^\infty$ by setting $\tau((b(f))=\hat{f}(0)$  \\
$(= (2 \pi)^{-d} \int f(x) dx)$, and easily verify (\cite{Fol}, page 
36) that it is a well-defined faithful trace on
 $\cla^\infty$. It is natural to consider $\clh=L^2(\cla^\infty, 
\tau)$ and represent $\cla$ in $\clb(\clh)$ by left
 multiplication. From the definition of $\tau$, it is clear that the 
map $C_c^\infty (\IR^{2d}) \ni f \mapsto b(f) \in
 \cla^\infty \subseteq \clh$ extends to a unitary isomorphism from 
$L^2(\IR^{2d})$ onto $\clh$ and in the sequel we shall often identify
 the two.

There is a canonical $2d$-paramater group of automorphism of $\cla$ 
given  by $\varphi_\al(b(f))=b(f_\al)$, where
 $\hat{f_\al}(x)=e^{i \al.x}\hat{f}(x),~f \in C_c^\infty 
(\IR^{2d}),~\al \in \IR^{2d}.$
Clearly, for any fixed $b(f)
\in \cla^\infty$, $\al \mapsto \varphi_\al(b(f))$ is smooth, and on 
differentiating
this map at $\al=0$, we get the
 canonical derivations $\dlt_j, j=1,2,\ldots, 2d$ as 
$\dlt_j(b(f))=b(\partial_j(f))$ for $f \in C_c^\infty (\IR^{2d}).$
 We shall not notationally distinguish between the derivation 
$\dlt_j$ on $\cla^\infty$ and its extension to $\clh$, and
  continue to denote by $i \dlt_j$ both the derivation on 
$\ast$-algebra $\cla^\infty$ and the associated
self-adjoint operator in
 $\clh$.

Let us now go back to the classical case. As a Riemannian manifold, 
$\IR^{2d}$ does not posses too many interesting
 features; it is a flat manifold and thus there is no nontrivial 
curvature form.  Instead, we shall be interested in
 obtaining the volume form from the operator-theoretic data 
associated with the $2d$-dimensional Laplacian $\Dlt$.
 Let $\clt_t=
e^{\frac {t}{2} \Dlt}$ be the contractive $C_0$-semigroup generated 
by $\Dlt$, called the heat semigroup on $\IR^{2d}$. Unlike
 compact manifolds, $\Dlt$ has only absolutely continuous spectrum. 
But for any $f \in C_c^\infty (\IR^{2d})$ and
 $\epsilon >0$, $M_f(-\Dlt +\epsilon)^{-d}$ has discrete spectrum. 
Furthermore, we have the following :
\bthm
\label{tr-vol}
$M_f \clt_t$ is trace-class  and $Tr(M_f \clt_t) = t^{-d} \int f(x) 
dx.$ Thus, in particular, $v(f)\equiv \int f(x)dx=
 t^d  Tr(M_f \clt_t).$
\ethm
{\it Proof :-}\\
We have $Tr(M_f \clt_t)=Tr(\clf M_f \clf^{-1} M_{e^{-\frac {t}{2} 
\sum x_j^2}})$, and $ \clf M_f \clf^{-1}
 M_{e^{-\frac {t}{2}\sum x_j^2}}$ is an integral
 operator with the kernel $k_t(x, y)= \hat{f}(x-y)e^{-\frac {t}{2} 
\sum y_j^2}.$ It is continuous in both arguments and $\int |k_t(x,x)|
dx < \infty$, we obtain by using a result in \cite{G-K}, (p. 114, 
ch.3) that $M_f \clt_t$ is trace class and $Tr(M_f \clt_t)=
 \int k_t(x,x)dx = {(2 \pi)}^d t^{-d} \hat{f}(0)= t^{-d} v(f).$
\qed

As in section 4, we get an alternative expression for the volume form 
$v$ in terms of the Dixmier trace.
\bthm
\label{dixvol}
For $\epsilon >0,$ $M_f (-\Dlt +\epsilon )^{-d}$ is of Dixmier trace 
class and its Dixmier trace is equal to
 $\pi^d v(f).$
\ethm
For convenience, we shall give the proof only in the case $d=1$. We 
need  following two lemmas.
\blmma
\label{L1}
If $f,g \in L^p(\IR^2)$ for some $p$ with $2 \leq p < \infty$, then 
$M_f \widetilde{M_g}$ is a compact operator in $L^2(\IR^2).$
 \elmma
{\it Proof :-}\\
It is a consequence of the Holder and Hausdorff-Young inequalities. 
We refer the reader to \cite{RS3}, volume III for a proof.
\qed

\blmma
\label{L2}
Let $S$ be a square in $\IR^2$ and $f$ be a smooth function with 
${\rm Supp}(f) \subseteq {\rm int}(S).$ Let $\Dlt_S$ denote
 the Laplacian on $S$ with the periodic boundary condition. Then 
$Tr_\omega (M_f (-\Dlt_S + \epsilon)^{-1})= \pi \int f(x) dx.$
 \elmma
{\it Proof :-}\\
This follows from  \cite{Land} by identifying $S$ with the 
two-dimensional torus in the natural manner.
\qed \\
{\it Proof of the theorem :-}\\
Note that for $g \in \cld(\Dlt) \subseteq L^2(\IR^2),$ we have $fg 
\in \cld(\Dlt_S)$
and $(\Dlt_S M_f -M_f \Dlt)(g)=
 (\Dlt M_f -M_f \Dlt)(g)=B g,$ where $B=-M_{\Dlt f} +2 i \sum_{j=1}^2 
M_{\partial_j(f) }
\circ \partial_j.$
 From this follows the identity
\bea
\label{**}
\lefteqn{M_f(-\Dlt+\epsilon)^{-1}-(-\Dlt_S +\epsilon)^{-1} 
M_f}{\nonumber}\\
&=& (-\Dlt_S +\epsilon)^{-1}B(-\Dlt+\epsilon)^{-1}.
\eea
 Now, from
 the Lemma \ref{L1}, it follows that $B(-\Dlt +\epsilon)^{-1}$ is 
compact, and since $(-\Dlt_S +\epsilon)^{-1}$ is of
 Dixmier trace class (by the Lemma \ref{L2}), we have that the right 
hand side of (\ref{**}) is of Dixmier trace class
 with the Dixmier trace $=0$. The theorem follows from the genaral 
fact that $Tr_\omega (xy)=Tr_\omega(yx)$, if $y$ is
 of Dixmier trace class and $x$ is bounded (see \cite{Con}).
\qed

 Similar  computation can be done for the non-commutative case. The 
Lindbladian  $\cll_0$ generated by the
canonical derivation $\dlt_j $ on $\cla$ is given by
\bea
\label {30}
\cll_0 (a(f)) = \frac{1}{2} a ( \Dlt f) ,~f \in C_c^\infty ( R^{2d} )
\eea
Since in $L^2(R^{2d} ) $ $\frac {1} {2} \Dlt$ has a natural 
selfadjoint extension   ( which we continue to express by the same 
symbol),
$\cll_0$ also has an extension as a negative selfadjoint operator in 
$\clh \cong L^2 ( R^{2d} ) $, and we define the
 heat semigroup for this case as $\clt_t = e^{t \cll_0} $. By analogy 
we can define the volume form on $\cla^\infty $ by setting
$ v(a(f)) = \lim_{t \raro 0+}   {t^d}  Tr (a(f) \clt_t ) $. Then we 
have

\bthm
$ v(a(f))= \int f dx $
\ethm
{\it Proof:--} \\
The kernel $\tilde {K}_t $ of the integral operator $a(f) \clt_t $ in 
$\clh$ is given as $ \tilde {K}_t (x,y)= \hat {f} (
\underline {x} - \underline {y} ) e^{-t|y|^2 /2} e^{ip(x,y)/2}$. As 
before we note that $K_t$ is continuous in $R^{2d} $
and $\tilde {K}_t (x,x) = k_t(x,x)= \hat {f} (0) e^{-t|x|^2 /2 }.$ 
Using \cite {G-K} we get the required result. \qed   \\

{\bf Remark}:--  (i) Note that in the theorem \ref{dixvol}, 
$Tr_{\omega} ( M_f (-\Dlt + \epsilon )^{-d}) = \pi^d v(f) $
which is independent of $\epsilon >0$. This could also have been 
arrived at directly as in section  4 for the algebra
$\ath$ once we have observed in the proof of the theorem that $ 
Tr_{\omega} M_f (\Dlt - \epsilon)^{-1}
= Tr_{\omega} M_f ( \Dlt_S -
\epsilon )^{-1} $.

We want to end this section with a brief discussion on the stochastic 
dilation of the heat semigroups on the spaces
 considered. For the classical (or commutative) $C^*$-algebra of $C_0 
(R^{2d}) $ the stochastic process
associated with the heat semigroup is the well known standard 
Brownian motion. For the non-commutative $C^*$-algebra  $\cla$
we first realize it in $\clb (L^2 (R^{d} ) )$ by the 
Stone-von-Neumann theorem on the representation of the Weyl relations
\cite{Fol}
\bea
\label {40}
  (U_{ {\al}} f)( {x})=f (  {x} +  {\bta} ){\nonumber}\\
(V_{ {\bta}}f )(  {x}) = e^{i {\al}.{x}} f
( {x}).
\eea
Let $ q_j,p_j ( j=1,2 \ldots d )$ be the generators of $ V_{ \bta} $ 
and $ U_{ \al}$
respectively, in fact they are the position and momentum operators in 
the above Schrodinger representation.
For simplicity of writing we shall restrict ourselves to the case 
$d=1$, and consider the q.s.d.e in $ L^2 (R)
\otimes \Gamma (L^2(R_+,C^2)):$
\bea
\label {41}
  dX_t= X_t [ -i p ~dw_1  (t) - \frac{1}{2}  p^2 dt   -i q~dw_1  (t) 
- \frac{1}{2}
 q^2 dt ] ,X_0=I
\eea
where $w_1,w_2$ are independent standard Brownian motions as in 
section 2. The following theorem
summarizes the results.
\bthm
(i) The q.s.d.e (\ref {41}) has a unique unitary solution. \\
(ii) If we set $j_t(x)= X_t (x \otimes I_t) X_t^* $ then $ j_t$ 
satisfies the  q.s.d.e :
$$ dj_t(x)  = j_t ( -i [p,x]) dw_1(t) +  j_t(-i [q,x])dw_2(t) 
+j_t(\cll (x) )dt $$ for all $ x \in \cla^\infty $
and $E j_t (x)= e^{t \cll } (x) $ for all $x \in \cla $
\ethm
{\it Proof:--} \\
Consider the q.s.d.e in $ \Gamma (L^2(R_+))$ for each $ \lmd \in R $ 
for a.a $w_1$,
$$ dW_t^{( \lmd)} = W_t^{(\lmd)} ( -i( \lmd +w_1 (t) )dw_2(t) - \frac 
{1} {2} (\lmd + w_1(t))^2 dt ), W_0^{(\lmd)}= I.$$
It is clear from \cite {KRP} that $W_t^{(\lmd)} = exp(-i \int_0^t 
(\lmd  + w_1(s) )dw_2 (s))$  which is unitary in
$\Gamma  (L^2(R_+))$ for fixed $\lmd$ and $w_1$. Next we set $W_t = 
\int_R E^q (d \lmd) \otimes W_t^{(\lmd)} $
which can be easily seen to be unitary in $ L^2 (R) \otimes \Gamma 
(L^2(R_+))$ for fixed $w_1$, where $E^q$ is the spectral
measure of the self adjoint operator $q$ in $L^2(R)$. Writing  $X_t = 
W_t e^{-ipw_1(t)} $ it  is clear that $ X_t $
is unitary in $L^2 (R) \otimes \Gamma (L^2(R_+,C^2))$. A simple 
calculation using Ito calculus shows that
 $X_t$ indeed satisfies equation \ref {41}.   \\
The part two follows from the observation that for fixed $w_1$ and 
$w_2$,  $X_t^*$ and $b(f) \otimes I_{\Gamma}$
with $ f \in C_c^{\infty} (R^2) $ maps $ \cls (R) \otimes \Gamma 
(L^2(R_+ , C^2 ))$ into itself. It is
also easy to see that  $j_t(x) = X_t x X_t^* = e^{-iqw_2(t)}  
e^{-ipw_1 (t)}  x e^{ipw_1(t)} e^{iqw_2(t)}
=\phi_{(-w_1(t),-w_2(t))}.$ \qed

\noindent {\it Acknowledgement :}\\
P. S. Chakraborty and D. Goswami acknowledge the support from the 
National Board of Higher Mathematics,
 India and K. B. Sinha acknowledges the support from the Jawaharlal 
Nehru Centre for Advanced
 Scientific Research, Bangalore, India.\\

\vspace{1in}

Stat-Math Division, Indian Statistical Institute, \\
7, S. J. S. S. Marg, New Delhi-110016, India.\\
e-mail : parthac@isical.ac.in,\\
goswamid@wiener.iam.uni-bonn.de,\\
kbs@isid.ac.in

\end{document}